\def\U{{\cal U}}
\def\V{{\cal V}}
\def\diam{\mathop{\rm diam}}
\def\conv{\mathop{\rm conv}}

\centerline{\bf Buildings have finite asymptotic dimension}
\medskip
\centerline{Jan Dymara\footnote{$^1$}{Instytut Matematyczny, Uniwersytet Wroc\l awski, pl. Grunwaldzki 2/4, 50-384 Wroc\l aw, 
Poland, \smallbreak{\tt dymara@math.uni.wroc.pl}. Partially supported by MNiSW grant N201 012 32/0718 and by DFG grant \smallbreak 436 POL 113/110/0-1.} 
and Thomas Schick\footnote{$^2$}{Mathem.~Institut, Georg-August-Universit\"at
G\"ottingen, Bunsenstr.~3, 37073 G\"ottingen, Germany,
\smallbreak{\tt schick@uni-math.gwdg.de}. Partially supported by DFG grant 436 POL 113/110/0-1.}
}
\bigskip
The goal of this note is to prove the following theorem, generalizing some
results of Matsnev [Theorem 3.22, Mat].
\medskip
\bf Theorem 1.\sl\par
The asymptotic dimension of any building is finite and equal
to the asymptotic dimension of an apartment in that building.
\rm\par
\medskip
Generally we use definitions and notation as in [D]. In particular,
$(W,S)$ is a finitely generated Coxeter system, $C$ is a building 
with Weyl group $W$, $|C|$ is the Davis realization of $C$.
We will, however, confuse the Coxeter group and its abstract Coxeter
complex, denoting both by $W$; 
in particular, $|W|$ denotes the Davis complex. The $W$-valued 
distance in $C$ will be denoted $\delta_C$, 
while $\delta$ will be the gallery distance (i.e., $\delta=\ell\circ\delta_C$,
where $\ell(w)$ is the shortest length of a word in generators $S$ representing $w$).
Basic properties of minimal galleries in buildings can be found in [R], [G]. 
We fix a chamber $B\in C$ and define the $B$-based \it folding map \rm
as $\pi\colon C\to W$, $\pi(c)=\delta_C(B,c)$. We also use 
$\pi$ for the geometric realization $|C|\to|W|$ of this map.
The word `building' in the statement of Theorem 1 can be understood
either as the discrete metric space $(C,\delta)$, or as the $CAT(0)$
metric space $(|C|,d)$ (these spaces are quasi--isometric).
Neither thickness nor local finiteness of $C$  is assumed.

Recall that a metric space $X$ has asymptotic dimension $\le n$ if 
for any  $d>0$ there exist $n+1$ families 
$\U^0,\ldots,\U^n$ of subsets of $X$ such that: (1) $\bigcup_i\U^i$ is a uniformly bounded 
cover of $X$, and (2) for every $i$ any two sets $U,U'\in\U^i$ are \it $d$-disjoint\rm:
$d(U,U'):=\inf\{d(x,y)\mid x\in U, y\in U'\}\ge d$.
Basic properties of asymptotic dimension
can be found in [BD]. 
We only need the definition and the fact that
Coxeter groups have finite asymptotic dimension (cf. [Theorem B, DJ]). 

\medskip
Now we start the proof of Theorem 1. 
Let $n$ be the asymptotic dimension of the apartment $W$ of $C$. Because $W$
embeds isometrically into $C$, the asymptotic dimension of $C$ is at least
$n$. We have to prove that its asymptotic dimension is $\le n$. 

Fix therefore $d>0$.
Let $\U^0,\ldots,\U^n$ be 
uniformly bounded  families of $2d$-disjoint sets in $|W|$
such that $\bigcup_i\U^i$ is a cover of $|W|$. For  
 $U\in \U^i$ let $N_d(U)=\{x\in |W|\mid d(x,U)<d\}$,
and let $\conv{N_d(U)}$ be  the convex hull of this set in 
the $CAT(0)$ space $|W|$. Note that $\diam(\conv{N_d(U)})\le\diam(U)+2d$.
Let ${\cal C}_U$ be the set of path--connected components of 
$\pi^{-1}(\conv{N_d(U)})$. Put $\V_U=\{\pi^{-1}(U)\cap A\mid 
A\in{\cal C}_U\}$ and $\V^i=\bigcup_{U\in\U^i}\V_U$.
It is obvious that $\bigcup_i\V^i$ is a cover of $|C|$.

\medskip
{\bf{Claim}}: $\V^i$ are uniformly bounded families
of $d$-disjoint sets. 

Evidently, the claim implies Theorem 1.
To prove the claim, we first establish $d$-disjointness. Let $x\in V\in\V_U$,
$x'\in V'\in\V_{U'}$ with $V\ne V'$ and $U,U'\in\U^i$.
There are two cases: $U\ne U'$ and $U=U'$. 
In the first case $d(x,x')\ge d(\pi(x), \pi(x'))$. But $\pi(x)\in U$,
$\pi(x')\in U'$ and $d(U,U')\ge 2d$.  
Therefore $d(V,V')\ge 2d>d$. 
For the second case, assume that $U=U'$. 
Suppose also that $V=\pi^{-1}(U)\cap A$ and $V'=\pi^{-1}(U)\cap A'$,
for some $A,A'\in{\cal C}_U$. Then the geodesic segment
$[x,x']$ is not entirely contained in $\pi^{-1}(\conv{N_d(U)})$
(otherwise $A$ and $A'$, hence $V$ and $V'$, coincide). 
Let $p\in[x,x']\setminus\pi^{-1}(\conv{N_d(U)})$;
we have $d(x,x')=d(x,p)+d(p,x')\ge
d(\pi(x),\pi(p))+d(\pi(p),\pi(x'))\ge d+d\ge d$.

It remains to check uniform boundedness. 
Let $V\in\V_U$; then $V\subseteq A$ for some path--connected component $A$ of 
$\pi^{-1}(\conv{N_d(U)})$. It is enough to find a uniform bound on the diameter 
of $A$. We would like to make $A$ gallery--connected. Since path--connected sets are usually
not gallery--connected, we perform an auxiliary thickening construction. 
Let $X$ be a subset of $|C|$ or of $|W|$. We put
$T(X)=\bigcup\{Res(p)\mid p\in X\}$, 
where $Res(p)=\bigcup\{|c|\mid p\in|c|\}$, which we call the {\it spherical
residue of $p$}.
Observe that $Res(p)$ is a geometric realization of a spherical building, and spherical buildings are 
gallery--connected; therefore, any two points in $A$ belong to chambers which can be connected
by a gallery in $T(A)$. Since $T(A)\subseteq\pi^{-1}\bigl(T(\conv{N_d(U)})\bigr)$, the set
$A$ is contained in a gallery--connected component of $\pi^{-1}\bigl(T(\conv{N_d(U)})\bigr)$.
Observe that,
uniformly in $U$,  
the diameter of $T(\conv{N_d(U)})$ is bounded by $R+2d+2\kappa$,
where  $R$ is the uniform bound on diameters of 
elements of the families $\U^i$, 
and $\kappa$ is the diameter of the realization of a chamber.
Since the distances $d$ and $\delta$ are quasi--isometric, it remains to prove the following lemma.
\smallbreak
\bf Lemma 1.\sl\par
For any $N>0$ there exists $M>0$ such that if $U$ is a subset of $W$ of 
$\delta$-diameter $\le N$, then any gallery--connected component $V$
of $\pi^{-1}(U)$ has $\delta$-diameter $\le M$.
\rm \smallskip
Lemma 1 follows from the next two lemmas.
\smallbreak
\bf Lemma 2.\sl\par
Let $W$ be a Coxeter group.
For any $R>0$ there exists $D=D(R)$ such that for any subset $U\subset |W|$ of diameter
$R$ satisfying  $d(|1|,U)>D$ there exists a codimension-one face of $|1|$ such that 
the wall containing that face separates $|1|$ from $U$.
\rm
\smallbreak
Proof. For $s\in S$, let $M_s$ be the wall containing the $s$-face of $|1|$,
and let $M_s^-$ be the open half-space with boundary $M_s$ 
that does not intersect $|1|$.
Let $r$ be greater than twice the diameter of a chamber, and let $b$ be the barycentre of $|1|$. 

We denote with $B_r(b)$ the open distance ball around $b$. We claim that 
$\{M_s^-\setminus B_r(b)\mid s\in S\}$ is an open cover of 
$|W|\setminus B_r(b)$.
If a point $q\in |W|\setminus B_r(b)$ is an interior point of a chamber, 
then we can consider a minimal gallery from $1$ to that chamber; this gallery starts by
crossing some wall $M_s$, and then $q\in M_s^-$. In general, we apply the above argument to the 
chamber $c$ in $Res(q)$ that is $\delta$-closest to $|1|$. 
Any other chamber $c'$ in $Res(q)$ can be connected
to $1$ by a minimal gallery passing through $c$, so that if a wall $M_s$ separates 
$1$ from  $c$ then it also separates $1$ from $c'$. Therefore $M_s$ separates $b$ from all 
points in the interior of $Res(q)$, in particular from $q$. 

Intersecting open sets from the family $\{M_s^-\setminus B_r(b)\mid s\in S\}$
with the distance sphere $S_r(b)$, we obtain an open cover of
this compact metric space. Let $\epsilon$ be the Lebesgue number of this
cover. Put $D=Rr/\epsilon$. We claim that 
the Lebesgue number of the cover $\{M_s^-\setminus B_{D}(b)\mid s\in S\}$ of $|W|\setminus B_{D}(b)$
is at least $R$. Indeed, let $U'$ be a subset of  $|W|\setminus B_{D}(b)$ of diameter $\le R$.
Because $W$ is a $CAT(0)$-space, the set $U$ of intersection points of $S_r(b)$
and geodesic intervals from $b$ to points in $U'$
has diameter $\le \epsilon$. Thus $U\subset M_s^-$ for some $s\in S$. Now if a point $q\in U'$ did 
not belong to $M_s^-$, then $S_r(b)\cap[b,q]$ would not be in $M_s^-$ either (because
$b,q$ belong to the convex set $|W|\setminus M_s^-$), contradicting  $U\subset M_s^-$.
\hfill\hbox{QED(Lemma 2)}

 \medskip
For $X\subseteq C$ or $X\subseteq W$ we denote by $N(X)$ the
set of all chambers that belong to $X$ or 
are adjacent to a chamber that 
belongs to $X$. In other words, $N(X)$ is the 1-neighbourhood of $X$ in 
the $\delta$-metric.
Let $U$ be a subset of $W$ of $\delta$-diameter $\le N$.
There exists $R>0$ depending only on $N$ (and $W$) such that the $d$-diameter
of $|N(U)|$ is $\le R$.
Iterated application of Lemma 2 provides a minimal gallery $\gamma=(1,w_1,\ldots,w_k)$ 
such that the wall between $w_i$ and $w_{i+1}$ separates $w_i$ from $N(U)$ 
and $d(|w_k|,|N(U)|)\le D$. Note that this separation property implies that every chamber which
belongs to $U$ can be joint to $1$ by a minimal gallery which 
is a concatenation
$\gamma\delta$, i.e.~extending $\gamma$. 

\bf Lemma 3. \sl\par
Let $U$ and $\gamma$ be given as just described. Recall that $B$ is the
``base'' chamber in $C$, with $\pi(B)=1$. For any chamber $c$ meeting
$\pi^{-1}(U)$ there is a minimal gallery from $B$ to $c$ whose $\pi$-projection
extends $\gamma$.

For any gallery--connected component $V\subset C$ of $\pi^{-1}(U)$ there exists a chamber 
$e\in\pi^{-1}(w_k)$ such that any minimal gallery from $B$  
to a chamber in $V$ whose $\pi$-projection prolongs $\gamma$ passes
through $e$.
\rm
\smallskip
Proof.
For each chamber $c\in N(V)$ let $P(c)$ be the set of all minimal galleries 
from $b$ to $c$ that are of the form $\Gamma\Delta$, where
$\pi(\Gamma)=\gamma$. We first show that this set is not empty. Choose an
arbitrary minimal gallery from $B$ to $c$ with $\pi$-projection $\alpha$, and a
minimal gallery of the form $\gamma\delta$ from $1$ to $\pi(c)$. $\alpha$
and $\gamma\delta$ are minimal galleries with the same extremities, therefore
are related by a sequence of Tits moves (cf. [Chap 4, Proposition 5, Bourb], 
or [Theorem 2.11, Ron]). 
This sequence lifts to
a sequence of moves relating the original minimal gallery to a minimal
$\Gamma\Delta$ with $\pi$-projection $\gamma\delta$, as required.

Next we claim that if $\Gamma_1\Delta_1,\Gamma_2\Delta_2\in P(c)$, then 
$\Gamma_1=\Gamma_2$. In fact, $\pi(\Delta_1)$ and $\pi(\Delta_2)$  
are minimal galleries in $W$ with the same extremities, therefore 
again are related by a sequence of Tits moves. As above, this sequence lifts to
a sequence of moves relating $\Delta_1$ and $\Delta'_2$ (and keeping
extremities fixed), where $\Delta'_2$ has the same $\pi$-projection as 
$\Delta_2$. Now $\Gamma_1\Delta'_2$ and $\Gamma_2\Delta_2$ have the same extremities
and the same $\pi$-projections, so that they coincide.  

Therefore we can define $e(c)$ as $\Gamma\cap\pi^{-1}(w_k)$ for some
$\Gamma\Delta\in P(c)$, and then $e(c)$ does not depend on the choice of $\Gamma\Delta$.
We will now show that, for $c\in V$, $e(c)$ does not depend on $c$. Since $V$ is gallery connected, it is enough 
to check this independence for adjacent
$c,c'\in V$. There are three cases:
\item{1)} {$\delta(w_k,\pi(c))=\delta(w_k,\pi(c'))$.
Then there is a chamber $c''\in N(V)$ adjacent to both $c$ and $c'$ whose $\pi$-projection is closer
to $w_k$ than that of $c$. And for any $\Gamma\Delta\in P(c'')$ we have 
$\Gamma\Delta c\in P(c)$, $\Gamma\Delta c'\in P(c')$.}
\item{2)}{ $\delta(w_k,\pi(c))=\delta(w_k,\pi(c'))-1$. Then for any $\Gamma\Delta\in P(c)$ we have
$\Gamma\Delta c'\in P(c')$.}
\item{3)}{ $\delta(w_k,\pi(c))=\delta(w_k,\pi(c'))+1$. This is symmetric to 2), we just switch $c$ and $c'$.}

Finally, $e=e(c)$, where $c\in V$, does not depend on $c$ and is as claimed.\hfill\hbox{QED(Lemma 3)}
\medskip

Let $L$ is the maximal gallery distance 
between chambers of $d$-distance $\le D+\kappa$, where $\kappa$ is the
diameter of a realization of a chamber. As a result of Lemma 3, every chamber
in $V$ is at $\delta$-distance $<L+N$ from $e$, 
hence $V$ has $\delta$-diameter $<2L+2N$. This proves Lemma 1 and therefore Theorem 1.

\medskip
The result leaves as an open task to determine the precise asymptotic
dimension of a Coxeter complex $W$ as above. These are particularly nice
$CAT(0)$-complexes. In this context, it is plausible to expect that the
asymptotic dimension coincides with the classical (microscopic)
dimension for a suitable choice or modification of the complex $W$. 

This coincidence has been established in model examples 
(coming from slightly different contexts), in particular for
many homogeneous manifolds
in [Corollary 3.6, CG] and in [Theorem 12, BD2]. 
For simply connected Riemannian
manifolds with curvature bounded by $c<0$, the result has been proved by
Grave [Theorem 6.20, Gra]. For spaces with cocompact action of the isometry group which are
Gromov hyperbolic, a related result has been obtained by Buyalo and
Lebedeva [BL]. It is not known to us whether these results extend to non-positive
curvature, or whether one can find counterexamples in this wider
class. However, observe that it is
not even known that the asymptotic dimension of every finite dimensional
$CAT(0)$-space is
finite.

Nonetheless, we conjecture that the answer to the following question is yes,
at least in many good cases. Given this belief, it would be even more
interesting to find counterexamples.

{\bf Question}: Is it true that for the Coxeter complex $W$ as considered in
this note, the asymptotic dimension and the virtual cohomological dimension,
i.e.~the dimension of its Bestvina complex [Best] coincide? Note that the
latter is also given by the dimension of the $CAT(0)$-boundary (compare
[BM]).

It follows from [Corollary 4.11, Dr] that the asymptotic dimension
of $W$ is not smaller than its virtual cohomological dimension.

\bigskip
{\bf Acknowledgments}.

We thank Alexander Dranishnikov and 
Jacek \smash{\hbox{\'Swi\hskip2.8pt\lower6.1pt\hbox{`}\hskip-5.8pt a\-tkowski}}
for useful comments. This research was
supported by German science foundation
grant 436 POL 113/110/0-1, whose support for a visit of the first author to
Georg-August-Universit\"at G\"ottingen, and of the second author to Banach
conference center, B\c edlewo, is gratefully acknowledged. Moreover, we thank 
these
institutions for their hospitality.

\medbreak

\bf References\rm
\smallskip
\item{[BD]} Bell G., Dranishnikov A., Asymptotic dimension in Bedlewo, 2005. 
arXiv: math.GR/0507570

\item{[BD2]} Bell, G. C., Dranishnikov, A. N., A Hurewicz-type theorem for
asymptotic dimension and applications to geometric group theory.
Trans. Amer. Math. Soc.  358  (2006),  no. 11, 4749--4764  

\item{[Best]} Bestvina, M.,
The virtual cohomological dimension of Coxeter groups. Geometric group theory,
Vol. 1 (Sussex, 1991), 19--23, 
London Math. Soc. Lecture Note Ser., 181, Cambridge Univ. Press, Cambridge,
1993. 

\item{[BM]} Bestvina, M., Mess, G.,
The Boundary of Negatively Curved Groups.
    Journal of the American Mathematical Society, Vol. 4, (1991), 469-481.

\item{[Bourb]} Bourbaki N., ``Groupes et algebres de Lie, chapitres IV-VI" 
Hermann 1968.

\item{[BL]} Buyalo, Sergei; Lebedeva, Nina. Dimensions of locally and
asymptotically self-similar spaces. arXiv: math.MG/0601744

\item{[CG]} Carlsson G., Goldfarb B., On homological coherence of discrete 
groups, J. Algebra 276 (2004), 502--514.
 
\item{[D]} Davis M.W., 
Buildings are CAT(0), 
in ``Geometry and Cohomology in Group Theory", 
edited by P. Kropholler, G. Niblo, R. Stohr,
LMS Lecture Notes \bf 252\rm, 
Cambridge Univ. Press, 
1998, pp.108--123.

\item{[Dr]} Dranishnikov A., Cohomological approach to asymptotic dimension.
arXiv: math.MG/0608215

\item{[DJ]} Dranishnikov A., Januszkiewicz T., Every Coxeter group acts amenably on a compact space,
Proceedings of the 1999 Topology and Dynamics Conference (Salt Lake City, UT), Topology Proc. \bf 24\rm
(1999), Spring, 135--141.  

\item{[G]} Garrett P., Buildings and classical groups, Chapman \& Hall, 1997.

\item{[Gra]} Grave, Bernd. Coarse geometry and asymptotic dimension,
Mathematica Gottingensis (2006). arXiv: math.MG/0601744.

\item{[Gr]} Gromov, Mikhael., Asymptotic invariants of infinite groups.
Geometric group theory, Vol. 2 (Sussex, 1991), 1--295, London
Math. Soc. Lecture Note Ser., \bf 182\rm,  Cambridge Univ. Press, Cambridge, 1993. 

\item{[Mat]} Matsnev D., The Baum--Connes conjecture and group actions 
on affine buildings, PhD thesis, Pennsylvania State University, 2005.


\item{[Ron]} Ronan M.,
Lectures on buildings, Perspectives in Mathematics vol. 7,
Academic Press 1989. 

\bye